\documentclass[global,final]{svjour}
\usepackage{graphics}
\usepackage{times}
\usepackage{amsmath,amssymb,amsbsy}
\usepackage[h]{esvect}

\usepackage{cmap,graphicx,parskip,graphpap,wrapfig}
\usepackage[noend]{algorithmic}

\usepackage{thumbpdf}
\usepackage[rightcaption]{sidecap}

\usepackage{soul,array,calc,url,ragged2e,fixltx2e}
\usepackage[stretch=15,shrink=15,step=3]{microtype}
\usepackage{ellipsis}

\newcommand{\subparagraph}{\paragraph}
\usepackage{titlesec}
\titlespacing{\section} {0pt}{2ex}{1ex plus .2ex}
\titlespacing{\subsection} {0pt}{2ex}{1ex plus .2ex}
\titlespacing{\subsubsection}{0pt}{2ex}{1ex plus .2ex}
\titlespacing{\paragraph} {0pt}{1ex}{0.5em}
\usepackage{url}

\journalname{Foundations of Computational Mathematics}

\begin{document}

\title{Discrete Lie Advection of Differential Forms}

\author{P. Mullen\inst{1} \and A. McKenzie\inst{1}\and D. Pavlov\inst{1} \and L. Durant\inst{1} \and\\ Y. Tong\inst{2}\and E. Kanso\inst{3} \and J. E. Marsden\inst{1} \and M. Desbrun\inst{1}\\}

\institute{Computing + Mathematical Sciences\\ California Institute of Technology\\  Pasadena, CA 91125, USA\\
\and Computer Science \& Engineering\\ Michigan State University\\ East Lansing, MI 48824, USA\\
\and Aerospace and Mechanical Engineering\\ University of Southern California\\ Los Angeles, CA 90089, USA}
\date{\small \it To appear in J.FoCM, 2011}

\maketitle

\begin{abstract}
In this paper, we present a numerical technique for performing Lie advection of arbitrary differential forms. Leveraging advances in high-resolution finite volume methods for scalar hyperbolic conservation laws, we first discretize the interior product (also called~\emph{contraction}) through integrals over Eulerian approximations of extrusions. This, along with Cartan's homotopy formula and a discrete exterior derivative, can then be used to derive a discrete Lie derivative. The usefulness of this operator is demonstrated through the numerical advection of scalar fields and 1-forms on regular grids.
\end{abstract}

\subclass{35Q35, 51P05,  65M08.}
\keywords{Discrete contraction -- Discrete Lie derivative -- Discrete differential forms -- Finite volume methods -- Hyperbolic PDEs}

\section{Introduction}

Deeply-rooted assumptions about smoothness and differentiability of most
continuous laws of mechanics often clash with the inherently discrete
nature of computing on modern architectures. To overcome this difficulty, a vast number
of computational techniques have been proposed to \emph{discretize} differential equations,
and numerical analysis is used to prove properties such as stability, accuracy, and
convergence.  However, many key properties of a mechanical system are characterized by its
symmetries and invariants (\emph{e.g.}, momenta), and
preserving these features in the computational realm can be of
paramount importance~\cite{GeoIntBook}, independent of the order
of accuracy used in the computations.  To this end, \emph{geometrically-derived} techniques
have recently emerged as valuable alternatives to traditional, purely
numerical-analytic approaches. In particular, the use of
differential forms and their discretization as cochains has been
advocated in a number of applications such as
electromagnetism~\cite{BossavitBook,StToDeMa2008,Hiptmair02}, discrete
mechanics~\cite{WestMarsden}, and even fluids~\cite{ETKSD07}.
\medskip

In this paper we introduce a finite volume based technique for solving the discrete Lie advection equation, ubiquitous in most advection phenomena: \begin{equation} \frac{\partial \boldsymbol{\omega}}{\partial t} + \boldsymbol{\mathcal{L}_X \omega} = 0, \label{eq:LieCanonical} \end{equation}
where $\boldsymbol{\omega}$ is an arbitrary discrete differential $k$-form~\cite{Arnold:2006:FEEC,Bochev2006,DKT06} defined on a discrete manifold, and $\boldsymbol{X}$ is a discrete vector field living on this manifold. Our numerical approach stems from the observation, developed in this paper, that the computational treatment of discrete differential forms share striking similarities with finite volume techniques~\cite{LevequeBook} and scalar advection techniques used in level sets~\cite{SethianBook,OsherFedkiwBook}. Consequently, we present a discrete interior product (or \emph{contraction}) computed using any of the $k$-dimensional finite volume methods readily available, from which we derive a numerical approximation of the spatial Lie derivative $\boldsymbol{\mathcal{L}_X}$ using a combinatorial exterior derivative.

\subsection{Background on the Lie Derivative}

The notion of Lie derivative $\boldsymbol{{\cal L}_X}$ in Elie
Cartan's Exterior Calculus~\cite{Cartan1945} extends the usual
concept of the derivative of a function along a vector field
$\boldsymbol{X}$. Although a formal definition of this operator can
be made purely algebraically (see~\cite{AMR}, $\S$5.3), its nature is
better elucidated from a dynamical perspective~\cite{AMR} ($\S$5.4).
Consequently, the spatial Lie derivative (along with its closely
related time-dependent version) is an underlying element in all
areas of mechanics: for example, the rate of strain tensor in
elasticity and the vorticity advection equation in fluid dynamics
are both nicely described using Lie derivatives.
\medskip

A common context where a Lie derivative is used to
describe a physical evolution is in the advection of \emph{scalar
fields}: a scalar field $\boldsymbol{\rho}$ being advected
in a vector field $\boldsymbol{V}$ can be written as:
$\boldsymbol{ \partial \rho}/\boldsymbol{\partial t} +
\boldsymbol{\mathcal{L}_{V}} \boldsymbol{\rho}=0.$ %
The case of divergence-free vector fields
(\emph{i.e.}, $\boldsymbol{\nabla\! \cdot\! V}\! =\! 0$) has been the
subject of extensive investigation over the past several decades
leading to several numerical schemes for solving these types of
hyperbolic conservation laws in various applications (see,
\emph{e.g.},~\cite{Shi:02:FVWENO,Dupont2003,Titarev:04:FVWENO,Iske2004,Pullin,ETKSD07,MOF}).
Chief among them are the so-called finite volume methods~\cite{LevequeBook},
including upwind, ENO, WENO, and high-resolution techniques. Unlike finite
difference techniques based on point values (\emph{e.g.}, \cite{EngquistOsher,ENO,WENO}),
such methods often resort to the conservative form of the advection equation
and rely on cell averages and the integrated fluxes in between. The integral
nature of these finite volume techniques will be particularly suitable in our
context, as it matches the foundations behind discrete versions of exterior
calculus~\cite{Bochev2006,Arnold:2006:FEEC}.

\medskip

While finite volume schemes have been successfully used for over a
decade, they have been used almost solely to advect scalar
fields, be they functions or densities, or systems thereof (e.g., components of tensor fields). To the authors'
knowledge, Lie advection of inherently non-scalar entities such as differential forms has yet to benefit from these advances, as differential forms are \emph{not} Lie advected in the same manner as scalar fields.

\subsection{Emergence of Structure-Preserving Computations}

Concurrent to the development of high-resolution methods for scalar
advection, structure-preserving geometric computational methods have
emerged, gaining acceptance among engineers as well as
mathematicians~\cite{IMA_book}. Computational
electromagnetism~\cite{BossavitBook,StToDeMa2008}, mimetic (or
natural) discretizations~\cite{Nicolaides1997,Bochev2006}, and more
recently Discrete Exterior Calculus (DEC,~\cite{Hirani:2003,DKT06})
and Finite Element Exterior Calculus (FEEC,~\cite{Arnold:2006:FEEC})
have all proposed similar discrete structures that
discretely preserve vector calculus identities to obtain improved
numerics. In particular, the relevance of exterior calculus
(Cartan's calculus of differential forms~\cite{Cartan1945}) and
algebraic topology (see, for instance,~\cite{Munkres1984}) to
computations came to light.
\medskip

Exterior calculus is a concise formalism to express differential and
integral equations on smooth and curved spaces in a consistent
manner, while revealing the geometrical invariants at play. At its
root is the notion of differential forms, denoting antisymmetric
tensors of arbitrary order. As integration of differential forms is
an abstraction of the measurement process, this calculus of forms
provides an intrinsic, coordinate-free approach particularly
relevant to concisely describe a multitude of physical
models that make heavy use of line, surface and volume integrals~\cite{Burke,AMR,LovelockRund,Flanders,MoritaBook,Carroll,Frankel}.
Similarly, many physical measurements, such as fluxes, are performed as specific
local integrations over a small surface of the measuring instrument.
Pointwise evaluation of such quantities does not have
physical meaning; instead, one should manipulate those quantities
only as geometrically-meaningful entities integrated over
appropriate submanifolds---these entities and their geometric properties are
embodied in discrete differential forms.

\medskip

Algebraic topology, specifically the notion of chains and cochains (see,
\emph{e.g.},~\cite{Whitney1957,Munkres1984}, has been used to
provide a natural discretization of these differential forms and to
emulate exterior calculus on finite grids: a set of values on
vertices, edges, faces, and cells are proper discrete versions of
respectively pointwise functions, line integrals, surface integrals,
and volume integrals. This point of view is entirely compatible with
the treatment of volume integrals in finite volume methods, or
scalar functions in finite element methods~\cite{Bochev2006};
but it also involves the "edge elements" and "facet elements"
as introduced in E\&M as special $H_{\operatorname{div}}$ and
$H_{\operatorname{curl}}$ basis elements~\cite{Nedelec86}.
Equipped with such discrete forms of arbitrary degree,
Stokes' theorem connecting differentiation and integration
is automatically enforced if one thinks of
differentiation as the dual of the boundary operator---a
particularly simple operator on meshes. With these basic building
blocks, important structures and invariants of the continuous
setting directly carry over to the discrete world, culminating in a
discrete Hodge theory (see recent progress in~\cite{ArFaWi2010}).
As a consequence, such a discrete exterior
calculus has, as we have mentioned, already proven useful in many
areas such as electromagnetism~\cite{BossavitBook,StToDeMa2008},
fluid simulation~\cite{ETKSD07}, surface parameterization~\cite{GY03},
and remeshing of surfaces~\cite{TACSD07} to mention a few.

\smallskip

Despite this previous work, the contraction and Lie derivative of arbitrary discrete
forms---two important operators in exterior calculus---have received
very little attention, with a few exceptions.  The approach in~\cite{Bossavit:2003}
(that we will review in \S\ref{sec:DynamicLie}) is to exploit the duality between
the extrusion and contraction operators, resulting in an integral definition of the
interior product that fits the existing foundations. While a
discrete contraction was derived using linear ``Whitney'' elements,
no method to achieve low numerical diffusion and/or high resolution
was proposed. Furthermore, the Lie derivative was not discussed.
More recently Heumann and Hiptmair~\cite{Heumann08} leveraged this work to suggest
an approach similar to ours in a finite element framework for Lie advection of forms
of arbitrary degree, however only $0$-forms were analyzed.

\subsection{Contributions}%
In this paper we extend the discrete exterior calculus machinery by introducing
discretizations of contraction and Lie advection with low numerical diffusion. Our work
can also be seen as an extension of classical numerical techniques for hyperbolic
conservation laws to handle advection of arbitrary discrete differential forms. In
particular, we will show that our scheme in 3D is a generalization of finite volume
techniques where not only \emph{cell-averages} are used, but also \emph{face-}
and \emph{edge-averages}, as well as \emph{vertex values}.

\section{Mathematical Tools}
Before introducing our contribution, we briefly review the existing
mathematical tools we will need in order to derive a discrete Lie
advection: after discussing our setup, we describe the necessary
operators of Discrete Exterior Calculus, before briefly reviewing the
foundations of finite volume methods for advection. In this paper
continuous quantities and operators are distinguished from their
discrete counterparts through a {\bf bold} typeface.

\subsection{Discrete Setup}

\paragraph{Space Discretization.}%
Throughout the exposition of our approach, we assume a regular
Cartesian grid discretization of space. This grid forms an orientable
$3$-manifold cell complex $K = (V,E,F,C)$ with vertex set $V =
\{ v_i \}$, edge set $E = \{ e_{ij} \}$, as well as face set $F $
and cell set $C$. Each cell, face and edge is assigned an arbitrary
yet fixed intrinsic orientation, while vertices and cells always
have a positive orientation. By convention, if a particular edge
$e_{ij}$ is positively oriented then $e_{ji}$ refers to the same
edge with negative orientation, and similar rules apply for higher
dimensional mesh elements given even vs. odd permutations of their
vertex indexing.

\paragraph{Boundary Operators.} \label{sec:boundary}
Assuming that mesh elements in $K$ are enumerated with an
arbitrary (but fixed) indexing, the incidence matrices of $K$ then
define the boundary operators. For example, we let $\partial^1$
denote the $|V| \times |E|$ matrix with $(\partial^1)_{ve} = 1$
(resp., $-1$) if vertex $v$ is incident to edge $e$ and the edge
orientation points towards (resp., away from) $v$, and zero
otherwise. Similarly, $\partial^2$ denotes the $|E| \times |F|$
incidence matrix of edges to faces with $(\partial^1)_{ef} = 1$
(resp., $-1$) if edge $e$ is incident to face $f$ and their
orientations agree (resp., disagree), and zero otherwise. The
incidence matrix of faces to cells $\partial^3$ is defined in a
similar way. See~\cite{Munkres1984} for details.

\subsection{Calculus of Discrete Forms}%
\label{sec:DEC}%
Guided by Cartan's exterior calculus of differential forms on smooth
manifolds, DEC offers a calculus on discrete manifolds that
maintains the covariant nature of the quantities involved.

\paragraph{Chains and Cochains.}
At the core of this computational tool is the notion of
\emph{chains}, defined as a linear combination of mesh elements; a
$0$-chain is a weighted sum of vertices, a $1$-chain is a weighted
sum of edges, etc. Since each $k$-dimensional cell has a well-defined notion
of boundary (in fact its boundary is a chain itself; the boundary of
a face, for example, is the signed sum of its edges), the boundary
operator naturally extends to chains by linearity. A \emph{discrete
form} is simply defined as the dual of a chain, or
\emph{cochain}, a linear mapping that assigns each
chain a real number. Thus, a $0$-cochain (that we will abusively
call a $0$-form sometimes) amounts to one value per $0$-dimensional cell,
such that any $0$-chain can naturally pair with this cochain. More
generally, $k$-cochains are defined by one value per $k$-cell, and
they naturally pair with $k$-chains. The resulting pairing of a
$k$-cochain $\alpha^k$ and a $k$-chain $\sigma_k$ is the discrete
equivalent of the integration of a continuous $k$-form
$\boldsymbol{\alpha}^k$ over a $k$-dimensional submanifold
$\boldsymbol{\sigma}_k$: %
$$\int_{\boldsymbol{\sigma}_k} \boldsymbol{\alpha}^k \equiv \langle \alpha^k,\sigma_k \rangle.$$
While attractive from a computational perspective due to their conceptual simplicity and
elegance, the chain and cochain representations are also deeply rooted in a theoretical
framework defined by H. Whitney~\cite{Whitney1957}, who introduced the Whitney and deRham
maps that establish an isomorphism between the cohomology of simplicial cochains and the
cohomology of Lipschitz differential forms. With these theoretical foundations,
chains and cochains are used as basic building blocks for direct
discretizations of important geometric structures such as the deRham
complex through the introduction of two simple operators.

\paragraph{Discrete Exterior Derivative. }%
The differential $\mathbf{d}$ (called exterior derivative) is an existing exterior calculus operator that we will need in our construction of
a Lie derivative. The discrete derivative $\mathrm{d}$ is
constructed to satisfy Stokes' theorem, which elucidates the duality
between the exterior derivative and the boundary operator. In the
continuous sense, it is written
\begin{equation} \label{EQ:Stokes}
\int_{\boldsymbol{\sigma}} \mathbf{d} \boldsymbol{\alpha} = \int_{ \boldsymbol{ \partial \sigma } } \boldsymbol{ \alpha } .
\end{equation}

Consequently, if $\alpha$ is a discrete differential $k$-form, then
the ($k$+$1$)-form d$\alpha $ is defined on any ($k$+$1$)-chain $\sigma$ by
\begin{equation}
\left\langle \mathrm{d} \alpha , \sigma \right\rangle = \left\langle \alpha ,
\partial \sigma \right\rangle ,
\end{equation}
where $\partial \sigma $ is the ($k$-chain) boundary of $\sigma$, as
defined in \S\ref{sec:boundary}. Thus the discrete
differential $\mathrm{d}$, mapping $k$-forms to ($k$+$1$)-forms, is
given by the co-boundary operator, the transpose of the signed
incidence matrices of the complex $K$; $\; \mathrm{d}_0 =
(\partial^1)^T$ maps $0$-forms to $1$-forms, $\mathrm{d}_1 =
(\partial^2)^T$ maps $1$-forms to $2$-forms, and more generally in $n$D,
$\mathrm{d}_k=(\partial^{k+1})^T.$ In relation to standard 3D vector
calculus, this can be seen as $\mathbf{d}_0 \equiv \boldsymbol{\nabla }$, $\mathbf{d}_1 \equiv \boldsymbol{\nabla
\times}$, and $\mathbf{d}_2 \equiv \boldsymbol{\nabla \cdot}$.  The fact that the
boundary of a boundary is empty results in $\mathrm{dd} = 0$, which
in turn corresponds to the vector calculus facts that $\boldsymbol{\nabla \times
\nabla} = \boldsymbol{\nabla \cdot \nabla \times} = 0$. Notice that this operator
is defined purely combinatorially, and thus does \emph{not} need a
high-order definition, unlike the operators we will introduce later.

\subsection{Principles of Finite Volumes}
\label{sec:FV} %

Given the integral representation of discrete forms used in the previous section, a last numerical tool we will need is a method for computing solutions to advection problems in integral form.  Finite volume methods were developed for exactly this purpose, and while we now provide a brief overview of this general procedure for completeness, we refer the reader to~\cite{LevequeBook}
and references therein for further details and applications. One
approach of finite volume schemes is to advect a function $u(x)$ by a
velocity field $v(x)$ using a Reconstruct-Evolve-Average (REA)
approach. In one dimension, we can define the cell average of a
function $u(x)$ over cell $C_i$ with width $\Delta x$ as
$$\bar{u}_i = \frac{1}{\Delta x} \int_{C_i}
 u(x) \; dx \quad i = 1, 2, \ldots, N.$$%
Given $k$ adjacent cell averages, the method will reconstruct a function such that the average of $p(x)$ in
each of the $k$ cells is equal to the average of $u(x)$ in those
cells.  High-resolution methods attempt to build a reconstruction such that it has only high-order error terms in smooth regions, while lowering the order of the reconstruction in favor of avoiding oscillations in regions with discontinuities like shocks.  Such adaptation can be done through the use of slope limiters or by changing stencil sizes using essentially non-oscillatory (ENO) and related methods.  This reconstruction can then be evolved by the velocity field and averaged back onto the Eulerian grid.

\smallskip

Another variant of finite volume methods is one that computes fluxes through cell boundaries.  Employing Stokes' theorem, the REA approach can be implemented by computing only the integral of the reconstruction which is evolved through each face, and then differencing the incoming and outgoing integrated fluxes of each cell to determine its net change in density.  It is this flux differencing approach that will be most convenient for deriving our discrete contraction operator, due to the observation that the net flux through a face induced by evolving a function forward in a velocity field is equal to the flux through the face induced by evolving the face \emph{backwards} through the same velocity field.  This second interpretation of the integrated flux is the same as computing the integral of the function over an extrusion of the face in the velocity field, as will be seen in the next section, and therefore we may use any of the wide range of finite volume methods to approximate integrals over extruded faces.

\section{Discrete Interior Product and Discrete Lie Derivative}
In keeping with the foundations of Discrete Exterior Calculus, we present
the continuous interior product and Lie derivative operators in
their ``integral'' form, \emph{i.e.}, we present continuous
definitions of $\mathbf{i}_{\boldsymbol{X}} \boldsymbol{\omega}$ and $\boldsymbol{\mathcal{L}}_{\boldsymbol{X}} \boldsymbol{\omega}$
\emph{integrated} over infinitesimal submanifolds: these integral forms will be particularly amenable
to discretization via finite volume methods and DEC as we discussed earlier.

\subsection{Towards a Dynamic Definition of Lie Derivative}
\label{sec:DynamicLie}

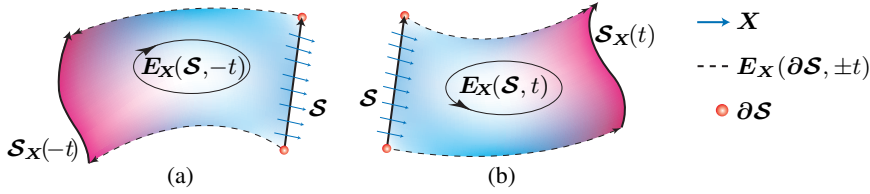
\begin{figure}[ht]\vspace*{-13mm}
\begin{picture}(200,100)
\put(0,0){\hspace*{5mm}\includegraphics[width=3.6in]{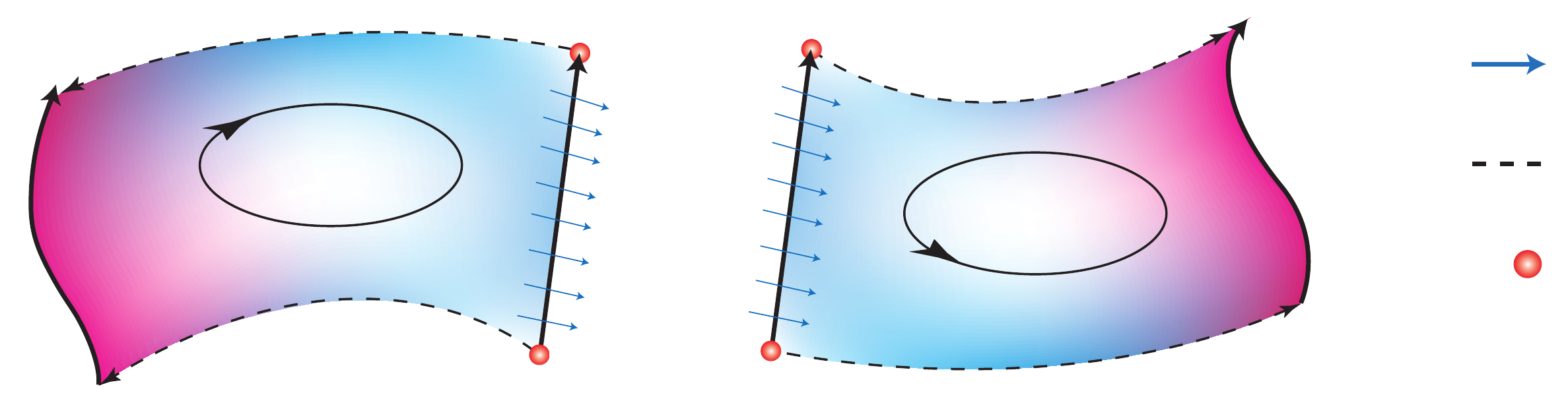}}
\put(273,53){$\boldsymbol{X}$}
\put(273,36){$\boldsymbol{E_X} (\boldsymbol{\mathcal{\partial S}},\pm t)$}
\put(273,19){$\boldsymbol{\partial \mathcal{S}}$}
\put(220,48){$\boldsymbol{\mathcal{S}_X}\!(t)$}
\put(169,29){$\boldsymbol{E\!_X}\!(\boldsymbol{\mathcal{S}},t)$}
\put(132,24){$\boldsymbol{\mathcal{S}}$}
\put(113,20){$\boldsymbol{\mathcal{S}}$}
\put(0,5){$\boldsymbol{\mathcal{S}_X}\!(\!-t\!)$}
\put(50,36){$\boldsymbol{E\!_X}\!(\boldsymbol{\mathcal{S}},\!-t)$}
\put(60,-5){(a)} \put(180,-5){(b)}
\end{picture}\vspace*{1mm}
\caption{\textbf{Geometric interpretation of the Lie derivative}
$\boldsymbol{\mathcal{L}}_{\boldsymbol{X}} \boldsymbol{\omega}$ of a differential form $\boldsymbol{\omega}$ in the
direction of vector field $\boldsymbol{X}$: (a) for a backwards advection in time of an edge $\mathcal{S}$ (referred to as upwind extrusion), and (b) for a forward advection of $\mathcal{S}$. Notice the orientation of the two extrusions are opposite, and depend on the direction of the velocity field.} \label{fig:continuousLie}
\end{figure}

\paragraph{Interior Product through Extrusion.}%
As pointed out by~\cite{Bossavit:2003}, the {\it extrusion} of
objects under the flow of a vector field can be used to give an
intuitive dynamic definition of the interior product. If
$\boldsymbol{\mathcal{M}}$ is an $n$-dimensional smooth manifold and $\boldsymbol{X} \in
\mathfrak{X}(\boldsymbol{\mathcal{M}})$ a smooth (tangent) vector field on the
manifold, let $\boldsymbol{\mathcal{S}}$ be a $k$-dimensional submanifold on
$\boldsymbol{\mathcal{M}}$ with $k<n$. The flow $\boldsymbol{\varphi}$ of the vector field $\boldsymbol{X}$
is simply a function $\boldsymbol{\varphi}$ : $\boldsymbol{\mathcal{M}} \times \mathbb{R}
\rightarrow \boldsymbol{\mathcal{M}}$ consistent with the one-parameter (time)
group structure, that is, such that $\boldsymbol{\varphi}( \boldsymbol{\varphi}(\boldsymbol{\mathcal{S}},
t), s) = \boldsymbol{\varphi}(\boldsymbol{\mathcal{S}}, s+t)$ with $\boldsymbol{\varphi}(\boldsymbol{\mathcal{S}},0) =
\boldsymbol{\mathcal{S}}$ for all $s, t \in \mathbb{R}$. Now imagine that
$\boldsymbol{\mathcal{S}}$ is carried by this flow of $\boldsymbol{X}$ for a time $t$; we
denote the resultant ``flowed-out'' submanifold $\boldsymbol{\mathcal{S}_X}(t)$,
which is equivalent to the image of $\boldsymbol{\mathcal{S}}$ under the mapping
$\boldsymbol{\varphi}$, \emph{i.e.}, $\boldsymbol{\mathcal{S}_X}(t) \equiv \boldsymbol{\varphi}(\boldsymbol{\mathcal{S}},
t)$. The extrusion $\boldsymbol{E_X}(\boldsymbol{\mathcal{S}},t)$ is then the ($k$+$1$)-dimensional
submanifold formed by the advection of $\boldsymbol{\mathcal{S}}$ over the time
$t$ to its final position $\boldsymbol{\mathcal{S}_X}(t)$: it is the ``extruded''
(or ``swept out'') submanifold. This can be expressed formally as a
union of flowed-out manifolds,
$$\boldsymbol{E_X}(\boldsymbol{\mathcal{S}},t) = \bigcup_{\tau \in [0,t]} \boldsymbol{\mathcal{S}_X}(\tau)$$
where the orientation of $\boldsymbol{E_X}(\boldsymbol{\mathcal{S}},t)$ is defined such that
$$\partial \boldsymbol{E_X} = \boldsymbol{\mathcal{S}_X}(t) - \boldsymbol{\mathcal{S}} - \boldsymbol{E_X}(\partial \boldsymbol{\mathcal{S}}, t).$$
These geometric notions are visualized in
Figure~\ref{fig:continuousLie}, where the submanifold $\boldsymbol{\mathcal{S}}$
is presented as a $1$-dimensional curve, flowed out to form
$\boldsymbol{\mathcal{S}_X}(t)$, or alternatively, extruded to form $\boldsymbol{E_X
}(\boldsymbol{\mathcal{S}},t)$.

\smallskip

Using this setup, the interior product $\mathbf{i}_{\boldsymbol{X}}$ of a time-independent form $\boldsymbol{\omega}$ evaluated on $\boldsymbol{\mathcal{S}}$ can now be defined through one of its most crucial properties, \emph{i.e.}, as the instantaneous change of $\boldsymbol{\omega}$ evaluated on $\boldsymbol{E_X}(\boldsymbol{\mathcal{S}},t)$, or more formally,
\begin{equation} \label{EQ:DynamicContraction} \int_{\boldsymbol{\mathcal{S}}}
\mathbf{i}_{\boldsymbol{X}} \boldsymbol{\omega} = \frac{d}{dt} \bigg|_{t=0}
\int_{\boldsymbol{E_X}(\boldsymbol{\mathcal{S}},t)} \boldsymbol{\omega}. \end{equation}%
While this equation is coherent with the discrete spatial picture, for the discrete Lie advection we will also wish to integrate $\mathbf{i}_{\boldsymbol{X}} \boldsymbol{\omega}$ over a small time interval. Hence, by
taking the integral of both sides of
Eq.~\eqref{EQ:DynamicContraction} over the interval $[0,\Delta t]$,
the first fundamental theorem of calculus gives us
\begin{equation} \label{EQ:DynamInterProd} \int_0^{\Delta t} \left[
\int_{\boldsymbol{\mathcal{S}_X}(t)} \mathbf{i}_{\boldsymbol{X}} \boldsymbol{\omega} \right] d
 t = \int_{\boldsymbol{E_X}(\boldsymbol{\mathcal{S}},\Delta t)} \boldsymbol{\omega},
\end{equation} which will be used later on for the discretization of the time-integrated
interior product.

\paragraph{Algebraic and Flowed-out Lie Derivative.}%
Using a similar setup, we can formulate a definition of Lie
derivative based on the flowed-out submanifold $\boldsymbol{\mathcal{S}_X}(t)$.
Remember that the Lie derivative is a generalization of the
directional derivative to tensors, intuitively describing the change
of $\boldsymbol{\omega}$ in the direction of $\boldsymbol{X}$. In fact, the Lie derivative
$\boldsymbol{\mathcal{L}_X \omega}$ evaluated on $\boldsymbol{\mathcal{S}}$ is equivalent to
the instantaneous change of $\boldsymbol{\omega}$ evaluated on $\boldsymbol{\mathcal{S}_X}(t)$,
formally expressed by \begin{equation} \label{EQ:continuousLie}
\int_{\boldsymbol{\mathcal{S}}} \boldsymbol{\mathcal{L}_X \omega} = \frac{d}{dt} \bigg|_{t=0}
\int_{\boldsymbol{\mathcal{S}_X}(t)} \boldsymbol{\omega}, \end{equation} as a direct consequence
of the Lie derivative theorem~\cite{AMR}(Theorem 6.4.1). As
before, we can integrate Eq.~\eqref{EQ:continuousLie} over a small
time interval $[0, \Delta t]$, applying the Newton-Leibnitz formula
to find
\begin{equation} \label{EQ:DynamLieDeriv}\int_0^{\Delta t} \left[
\int_{\boldsymbol{\mathcal{S}_X}(t)} \boldsymbol{\mathcal{L}_X \omega }\right] d t =
\int_{\boldsymbol{\mathcal{S}_X}(\Delta t)} \boldsymbol{\omega} - \int_{\boldsymbol{\mathcal{S}}} \boldsymbol{\omega}.
\end{equation}%
Note that the formulation above, discretized using a
semi-Lagrangian method, has been used, \emph{e.g.},  by~\cite{ETKSD07} to advect
fluid vorticity; in that case the right hand side  of Eq.~\eqref{EQ:DynamLieDeriv}
was evaluated by looking at the circulation through the boundary of
the ``backtracked'' manifold. Rather than following their approach,
we revert to discretizing the dynamic definition of the interior
product in Eq.~\eqref{EQ:DynamInterProd} instead, and later
constructing the Lie derivative algebraically. The primary
motivation behind this modification is one of effective numerical
implementation: we can apply a dimension-by-dimension finite volume scheme to obtain an approximation of the interior product, while the alternative---computing
integrals of approximated $\boldsymbol{\omega}$ over $\boldsymbol{\mathcal{S}_X}(t)$ as required by
a discrete version of Eq.~\eqref{EQ:DynamLieDeriv}---is comparatively cumbersome.  Also, by building on top of standard finite volume schemes the solvers can leverage pre-existing code, such as CLAWPACK~\cite{clawpack}, without requiring modification.

\smallskip
We now show how the Lie derivative and the interior product are linked
through a simple algebraic relation known as Cartan's homotopy formula.
In particular, this derivation (using Figure~\ref{fig:continuousLie}
as a reference) requires repeated application of Stokes' theorem from
Eq.~\eqref{EQ:Stokes}.\allowdisplaybreaks
\begin{align}\label{EQ:proof1}
\lim_{\Delta t \rightarrow 0} \frac{1}{\Delta t} \int_0^{\Delta t}
\left[ \int_{\boldsymbol{\mathcal{S}_X}(t)}
\!\!\!\!\!\! \boldsymbol{\mathcal{L}_X \omega} \right] d t & =
\lim_{\Delta t \rightarrow 0} \frac{1}{\Delta t} \left[
\int_{\boldsymbol{\mathcal{S}_X}(\Delta t)} \!\!\!\!\!\!
\boldsymbol{\omega} - \int_{\boldsymbol{\mathcal{S}}} \boldsymbol{\omega}
\right]\\ \label{EQ:proof2} & = \lim_{\Delta t \rightarrow 0}
\frac{1}{\Delta t} \left[ \int_{\boldsymbol{E_X} (\boldsymbol{\mathcal{S}},\Delta t)}
\!\!\!\!\!\! \mathbf{d} \boldsymbol{\omega} + \int_{\boldsymbol{E_X} (\boldsymbol{\partial \mathcal{S}}, \Delta t)}
\!\!\!\!\!\! \boldsymbol{\omega} \right]\\ \label{EQ:proof3}& = \int_{\boldsymbol{\mathcal{S}}}
\mathbf{i}_{\boldsymbol{X}} \mathbf{d} \boldsymbol{\omega} +
\int_{\boldsymbol{\mathcal{\partial S}}} \mathbf{i}_{\boldsymbol{X}} \boldsymbol{\omega} \\
\label{EQ:proof4} & = \int_{\boldsymbol{\mathcal{S}}} \mathbf{i}_{\boldsymbol{X}} \mathbf{d}
\boldsymbol{\omega} + \int_{\boldsymbol{\mathcal{S}}} \mathbf{d} \mathbf{i}_{\boldsymbol{X}} \boldsymbol{\omega}.
\end{align}

The submanifolds $\boldsymbol{\mathcal{S}}$ and $\boldsymbol{\mathcal{S}_X}(\Delta t)$ form a portion of the boundary of $\boldsymbol{E_X} (\boldsymbol{\mathcal{S}}, \Delta t)$. Therefore, by Stokes', we can evaluate $\mathbf{d} \boldsymbol{\omega}$ on the extrusion and subtract off the other portions of $\boldsymbol{\partial E_X} (\boldsymbol{\mathcal{S}},\Delta
t)$ to obtain the desired quantity. This is how we proceed from Eq.~\eqref{EQ:proof1} to Eq.~\eqref{EQ:proof2} of the proof. The following line, Eq.~\eqref{EQ:proof3}, is obtained by applying the dynamic definition of the interior product given in Eq.~\eqref{EQ:DynamInterProd} to each of two terms, leading us to our final result in Eq.~\eqref{EQ:proof4} through one final application of Stokes' theorem. What we have obtained is the Lie derivative expressed algebraically in terms of the exterior derivative and interior product. Notice that Eq.~\eqref{EQ:proof4} is the integral form of the celebrated identity called Cartan's homotopy (or \emph{magic}) formula, most frequently written as
\begin{equation}\label{eq:MAGIC}
\boldsymbol{\mathcal{L}_X \omega} = \mathbf{i}_{\boldsymbol{X}} \mathbf{d} \boldsymbol{\omega} + \mathbf{d}
\mathbf{i}_{\boldsymbol{X}} \boldsymbol{\omega}.
\end{equation}%
By defining our discrete Lie derivative through this relation, we
ensure the algebraic definition holds true in the discrete sense by
construction. It also implies that the Lie derivative can be
directly defined through interior product and exterior derivative, without the need for its own discrete definition.

\paragraph{Upwinding the Extrusion}
We may rewrite the above notions using an ``upwinded'' extrusion (\emph{i.e.}, a cell extruded backwards in time)
as well (see Fig.~\ref{fig:continuousLie}a).  For example, Eq.~\eqref{EQ:DynamicContraction} can be
rewritten as
\begin{equation} \label{EQ:DynamicUpwindContraction}\int_{\boldsymbol{\mathcal{S}}}
\mathbf{i}_{\boldsymbol{X}} \boldsymbol{\omega} = -\frac{d}{dt} \bigg|_{t=0}
\int_{\boldsymbol{E_X}(\boldsymbol{\mathcal{S}},-t)} \boldsymbol{\omega}.\end{equation}
While this does not change the instantaneous value of the contraction, integrating Eq.~\eqref{EQ:DynamicUpwindContraction} over the time interval $[0,\Delta t]$ now gives us \begin{equation} \label{EQ:DynamUpwindInterProd}\int_0^{\Delta t} \left[
\int_{\boldsymbol{\mathcal{S}_X}(t)} \mathbf{i}_{\boldsymbol{X}} \boldsymbol{\omega} \right] d
 t = -\int_{\boldsymbol{E_X}(\boldsymbol{\mathcal{S}},-\Delta t)} \boldsymbol{\omega}.\end{equation}  Similar treatment for the remainder of the above can be done and Cartan's formula can be derived the same way, however by using these definitions in our following discretization we will obtain computations over upwinded regions equivalent to those computed by finite volume methods.

\subsection{Discrete Interior Product}
A discrete interior product is computed by exploiting the
principles of Eq.~\eqref{EQ:DynamInterProd} and applying the finite volume machinery. Given a discrete $k$-form $\alpha$
and a discrete vector field $X$, the interior product is
approximated by extruding backwards in time every ($k$-$1$)-dimensional cell
$\sigma$ of the computational domain to form a new $k$-dimensional cell $E_X
(\sigma,-\Delta t)$. Evaluating the integral of $\alpha$ over the
extrusion and assigning the resulting value to the original cell
$\sigma$ yields the mapping $\left\langle \mathrm{i}_X \alpha,
\sigma \right\rangle$ integrated over a time step $\Delta t$. This
procedure, once applied to all ($k$-$1$)-dimensional cells, gives the desired discrete ($k$-$1$)-form $\mathrm{i}_X \alpha$.

\begin{figure}[h]\vspace*{-7mm}
\begin{picture}(170,100)
\put(35,0){\includegraphics[width=0.8\textwidth]{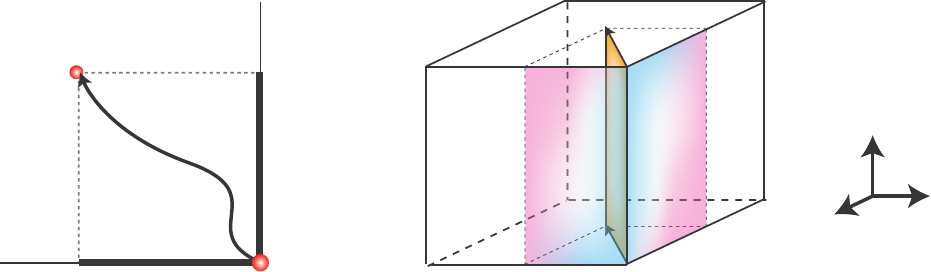}}
\put(20,10){(a)}
\put(140,10){(b)}
\put(267,8){$z$}
\put(292,25){$x$}
\put(285,35){$y$}
\end{picture}
\caption{\textbf{Approximating Extrusions}: In the discrete setting, the extrusion of a $(k$-$1)$-dimensional
manifold ($k$=$1$ on left, $2$ on right) is
approximated by projecting the Lagrangian advection of the manifold
into $\binom{n}{k}$ separate $k$-dimensional components. \vspace*{-1mm}}
\label{fig:projection}
\end{figure}
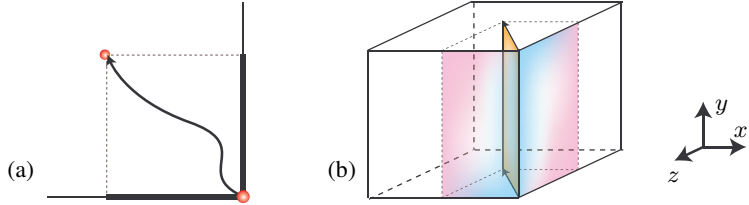

\smallskip

\paragraph{K-dimensional Splitting.}
One option for computing this integral would be to do an $n$-dimensional reconstruction of $\alpha$, perform a Lagrangian advection of the cell $\sigma$ to determine $E_X(\sigma, -\Delta t)$, and then algebraically or numerically computing the integral of the reconstructed $\alpha$ over this extrusion.  In fact, this is the idea behind the approach suggested in~\cite{Heumann08}.  However, with the exception of when $k \! = \! n$, such an approach does not allow us to directly leverage finite volume methods, as performing an $n$-dimensional reconstruction of a form given only integrals over $k$-dimensional submanifolds would require a more general finite element framework. For simplicity and ease of implementation we avoid this generalization and instead resort to projecting the extrusion onto the grid-aligned $k$-dimensional subspaces and then applying a $k$-dimensional finite volume method to each of the $\binom{n}{k}$ projections.  The integrals over the extrusion of $\sigma$ from each dimension are then summed.  Again, note that in the special case of $k \! = \! n$ no projection is required and we are left exactly with an $n$-dimensional finite volume scheme.  We have found that this splitting combined with a high-resolution finite volume method, despite imposing at most first order accuracy, can still give high quality results with low numerical diffusion, while being able to leverage existing finite volume solvers without modification.  However, if truly higher order is required then a full-blown finite element method would most likely be required~\cite{Heumann08}.

\paragraph{Finite Volume Evaluation.}
As hinted at in Section~\ref{sec:FV}, we notice that the time integral of the flux of a density field being advected through a submanifold $\sigma$ is equivalent to the integral of the density field over the backwards extrusion of $\sigma$ over the same amount of time.  In fact, some finite volume methods are derived using this interpretation, doing a reconstruction of the density field, approximating the extrusion, and integrating the reconstruction over this.  However, many others are explained by computing a numerical flux per face, and then multiplying by the time step $\Delta t$: this is still an approximation of the integral over the extrusion, taking the reconstruction to be a constant (the numerical flux divided by $v$) and the extrusion having length $v \Delta t$.  Indeed the right hand side of Eq.~\eqref{EQ:DynamicContraction} can be seen as analogous to the numerical flux, after which Eq.~\eqref{EQ:DynamInterProd} becomes the relationship between integrating the flux over time and the form over the extrusion.  Hence we may use any of the finite volume methods for $k$-dimensional density advection problems when computing the contraction of a $k$-form.  The only difference here is that rather than applying Stokes theorem and summing the contributions back to the original $k$-cell (which will be done by the discrete exterior derivative in the $\mathbf{d} \mathbf{i}_{\boldsymbol{X}} \boldsymbol{\omega}$ term of the Lie derivative), the contraction simply stores the values on the ($k$-$1$)-cells, without the final sum.

\subsection{Discrete Lie Advection}
We now have all the ingredients to introduce a discrete Lie
advection. Given a $k$-form $\alpha$, we compute the ($k$+$1$)-form
$\mathrm{d} \alpha$ by applying the transpose of the incidence
matrix $\partial^{k+1}$ to $\alpha$ as detailed in \S\ref{sec:DEC}. We
then compute the $k$-form $\mathrm{i}_X(\mathrm{d}\alpha)$, and the
($k$-$1$)-form $\mathrm{i}_X\alpha$. By applying
$\mathrm{d}$ to the latter form and summing the resulting
$k$-form with the other interior product, we finally get an approximation
of Cartan's homotopy formula of the Lie derivative. An explicit example of
this will be given in the next section to better illustrate the process
and details.

\section{Applications and Results}
We now present a few direct applications of this discrete Lie
advection scheme.  In our tests we used upwinding one-dimensional WENO schemes for our contraction operator, splitting even the $k$-dimensional problems into multiple one-dimensional ones.  We found that when using high-resolution WENO schemes we could obtain quality results with little numerical smearing despite this dimensional splitting.

\paragraph{A Note on Vector Fields.} In this section we assume that vector fields are discretized by storing their flux (\emph{i.e.}, contraction with the volume form) on all the $(n\!-\!1)$-dimensional cells of a $n$D regular grid, much like the Marker-And-Cell ``staggered'' grid setup~\cite{MACgrids}.  Evaluation of the vector fields at lower dimensional cells is done through simple averaging of adjacent discrete fluxes. We pick this setup as it is one of the most commonly-used representations, but the vector fields can be given in arbitrary form with only minor implementation changes.

\subsection{Volume Forms and $0$-Forms}
Applying our approach to volume forms ($n$-forms in $n$ dimensions) we have $$\boldsymbol{\mathcal{L}_X \omega} = \mathbf{i}_{\boldsymbol{X}} \mathbf{d} \boldsymbol{\omega} + \mathbf{d}
\mathbf{i}_{\boldsymbol{X}} \boldsymbol{\omega} = \mathbf{d}
\mathbf{i}_{\boldsymbol{X}} \boldsymbol{\omega}.$$
Note that $\mathbf{i}_{\boldsymbol{X}} \boldsymbol{\omega}$ is the numerical flux computed by the chosen $n$-dimensional finite volume scheme while $\mathbf{d}$ will then just assign the appropriate sign of this flux to each cell's update, and hence we trivially arrive at the chosen finite volume scheme with no modification. Similarly, applying this approach to $0$-forms results in well-known finite difference advection schemes of scalar fields.  Indeed, we have in this case $$\boldsymbol{\mathcal{L}_X \omega} = \mathbf{i}_{\boldsymbol{X}} \mathbf{d} \boldsymbol{\omega} + \mathbf{d}
\mathbf{i}_{\boldsymbol{X}} \boldsymbol{\omega} = \mathbf{i}_{\boldsymbol{X}} \mathbf{d} \boldsymbol{\omega}$$ as the contraction of a $0$-form vanishes.  We are thus left with $\mathbf{d} \boldsymbol{\omega}$ computing standard finite differences of a node-based scalar field on edges, and $\mathbf{i}_{\boldsymbol{X}}$ then doing componentwise upwind integration of reconstructions of these derivatives.  Such techniques are common in scalar field advection, for example in the advection of level sets, and we refer the reader to~\cite{OsherFedkiwBook,SethianBook} and references therein for examples.

\subsection{Advecting a $1$-Form in 2D}%
The novelty of this approach comes when applied to $k$-forms in $n$ dimensions with $n > k > 0$.  We first demonstrate the simplest such application of our method by advecting a $1$-form by a static velocity field in 2D using the simple piecewise-constant upwinding finite volume advection.  To illustrate the general approach we will explicitly write out the algorithm for this case.  We will assume the velocity $\boldsymbol{X}$ is everywhere positive in both $x$ and $y$ components to simplify the upwinding, and $X^x$ and $X^y$ will be used to represent the integrated flux through vertical and horizontal edges respectively.
\begin{figure}
\begin{center}
\includegraphics[width=0.6\textwidth]{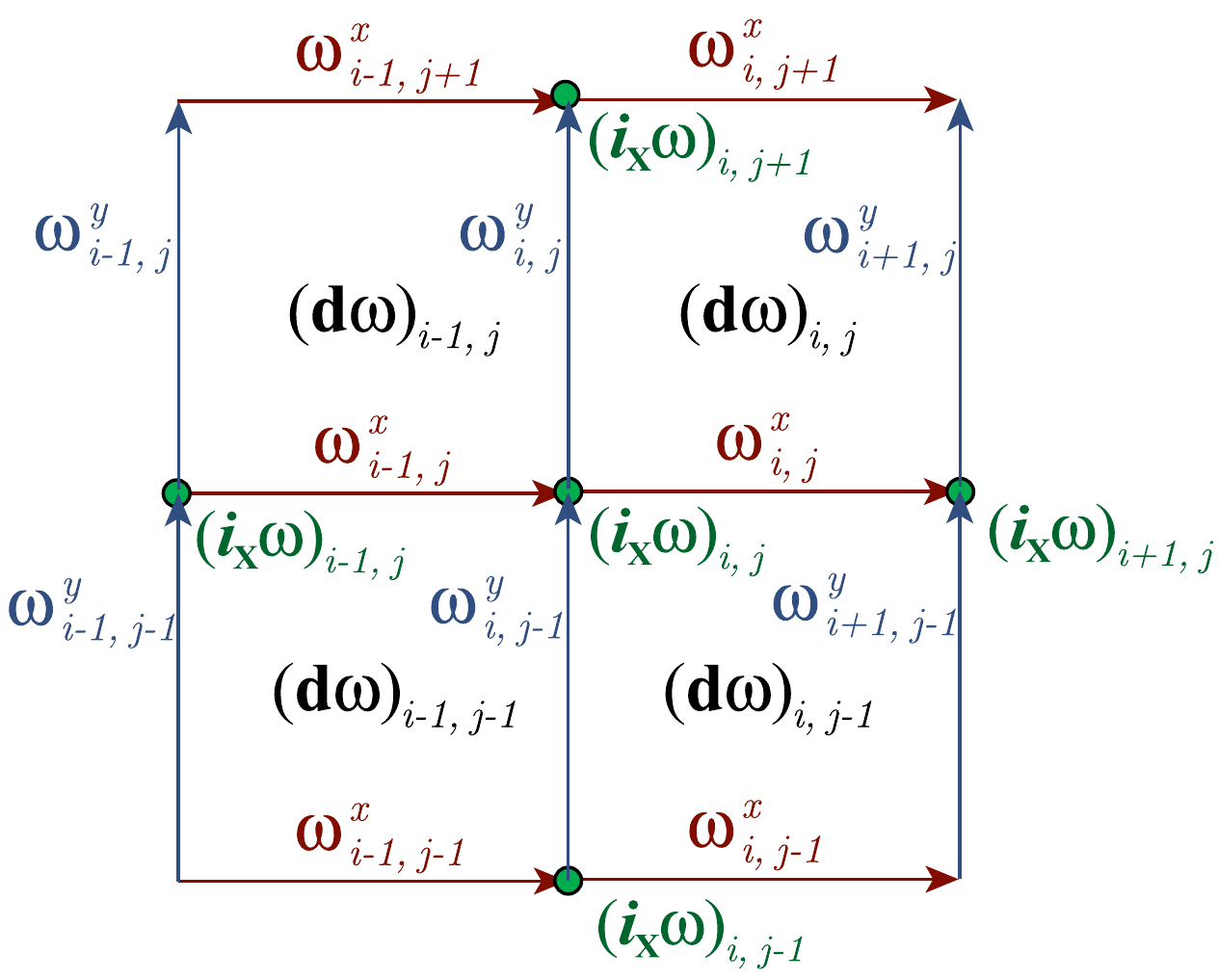}\vspace*{-1mm}
\end{center}
\caption{\textbf{Grid setup:} Indexing and location of the various quantities stored on different parts of the grid. Arrows indicate the orientation of the edges.}
\label{fig:gridsetup}
\end{figure}

\begin{figure}[h]\vspace*{-2mm}
\includegraphics[width=\textwidth]{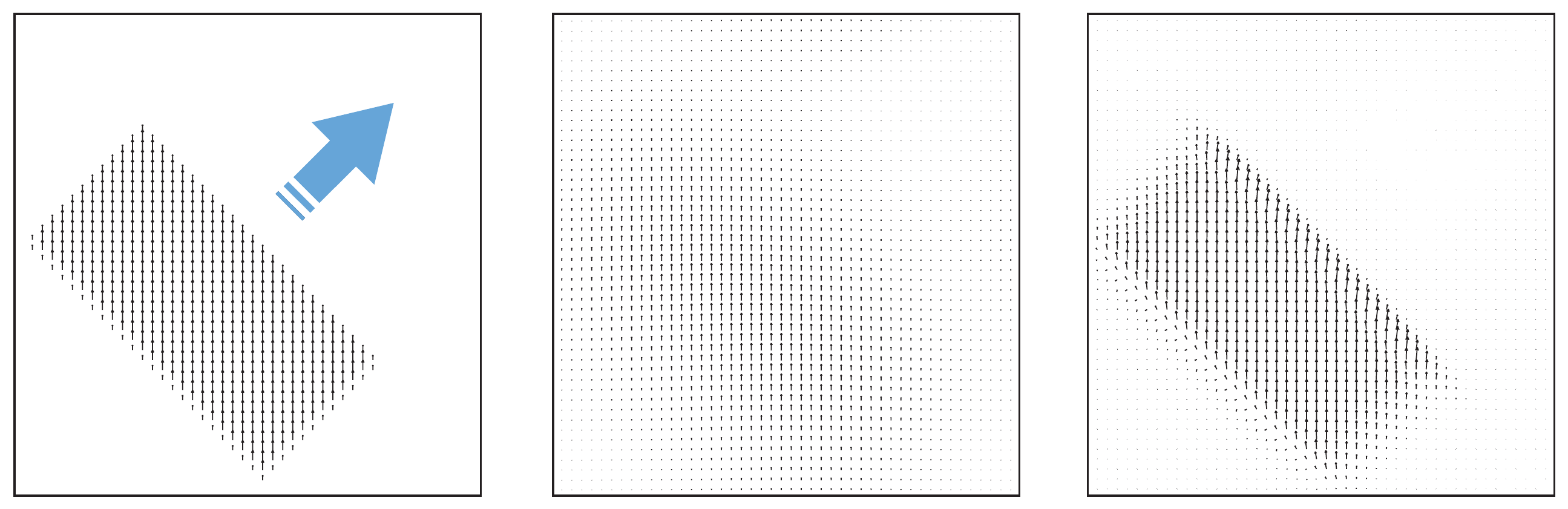}\vspace*{-1mm}
\hspace*{.14\textwidth}(a)\hspace*{.32\textwidth}(b)\hspace*{.32\textwidth}(c)\vspace*{-1mm}
\caption{\textbf{$\mathbf{1}$-Form advection:} (a) A piecewise-constant form ($dy$ within a rectangular shape, $0$ outside) is advected in a constant velocity field ($\mathbf{X}=(1,1)$, blue arrow) on a unit square periodic domain with a grid resolution of $48^2$ and a time step $dt\!=\!10^{-3}$. (b) Because the domain is periodic, the form should be advected back to its original position after $1s$ ($1000$ steps); however, our numerical method with a piecewise constant upwind finite volume scheme results in considerable smearing instead. (c) Using a high-resolution scheme (here, WENO-7) as the basic component of our form advection procedure significantly reduces smearing artifacts (same number of steps and step size).
\vspace*{-1mm}} \label{fig:compare0-7}
\end{figure}

Suppose we have a regular two-dimensional grid with square cells of size $h^2$, and with each horizontal edge oriented in the positive $x$ direction and each vertical edge oriented in the positive $y$ direction and numbered according to Fig~\ref{fig:gridsetup}.  A discrete $1$-form $\boldsymbol{\omega}$ is represented by its integral along each edge.  Due to the Cartesian nature of the grid, this implies that the $dx$ component of the form will be stored on horizontal edges and the $dy$ component will be stored on vertical edges, and we represent these scalars as $\boldsymbol{\omega}^x_{i,j}$ and $\boldsymbol{\omega}^y_{i,j}$ for the integrals along the $(i,j)$ horizontal and vertical edge respectively.  The discrete exterior derivative integrated over cell $(i,j)$, $(\mathbf{d} \boldsymbol{\omega})_{i,j}$, consists of the signed sum of $\boldsymbol{\omega}$ over cell $(i,j)$'s boundary edges, namely
$$(\mathbf{d} \boldsymbol{\omega})_{i,j} = \boldsymbol{\omega}^x_{i,j} + \boldsymbol{\omega}^y_{i+1,j} - \boldsymbol{\omega}^x_{i,j+1} - \boldsymbol{\omega}^y_{i,j}.$$
Using piecewise-constant upwind advection, and remembering the assumption of positivity of the components of $\boldsymbol{X}$,  we may now compute $\mathbf{i}_{\boldsymbol{X}} \mathbf{d} \boldsymbol{\omega}$ over a time interval $\Delta t$ for the horizontal and vertical edges $(i,j)$ as
\begin{equation}
\label{eq:ixdExample}
\begin{array}{rccl}
(\mathbf{i}_{\boldsymbol{X}} \mathbf{d} \boldsymbol{\omega})^x_{i,j} & = & - & \dfrac{\Delta t}{h^2} X^y_{i,j} (\mathbf{d} \boldsymbol{\omega})_{i,j-1}\\[1.5mm]
(\mathbf{i}_{\boldsymbol{X}} \mathbf{d} \boldsymbol{\omega})^y_{i,j} & = & & \dfrac{\Delta t}{h^2} X^x_{i,j} (\mathbf{d} \boldsymbol{\omega})_{i-1,j}.
 \end{array}
\end{equation}
Note the sign difference is due to the orientation of the extrusions, and would be different if the velocity field changed sign (see Fig.~\ref{fig:continuousLie}).  To compute the second half of Cartan's formula we must now compute $\mathbf{i}_{\boldsymbol{X}} \boldsymbol{\omega}$ at nodes, and then difference them along the edges.  Using dimension splitting, as well averaging the velocity field from edges to get values at nodes, we get for node $(i,j)$
\begin{equation}
\label{eq:ixExample}
(\mathbf{i}_{\boldsymbol{X}}\boldsymbol{\omega})_{i,j} = \frac{\Delta t}{2h^2} \left((X^x_{i,j} + X^x_{i,j-1}) \boldsymbol{\omega}^x_{i-1,j} + (X^y_{i,j} + X^y_{i-1,j}) \boldsymbol{\omega}^y_{i,j-1}\right).
\end{equation}
We may now trivially compute $\mathbf{d} \mathbf{i}_{\boldsymbol{X}} \boldsymbol{\omega}$ for edges as
\begin{align*}
(\mathbf{d} \mathbf{i}_{\boldsymbol{X}} \boldsymbol{\omega})^x_{i,j} & = (\mathbf{i}_{\boldsymbol{X}}\boldsymbol{\omega})_{i+1,j} - (\mathbf{i}_{\boldsymbol{X}}\boldsymbol{\omega})_{i,j}\\
(\mathbf{d} \mathbf{i}_{\boldsymbol{X}} \boldsymbol{\omega})^y_{i,j} & = (\mathbf{i}_{\boldsymbol{X}}\boldsymbol{\omega})_{i,j+1} - (\mathbf{i}_{\boldsymbol{X}}\boldsymbol{\omega})_{i,j}.
\end{align*}  Cartan's formula and the definition of Lie advection now lead us to obtain $$\Delta \boldsymbol{\omega} = - \int_0^{\Delta t}
\boldsymbol{\mathcal{L}_X \omega } d t$$ discretized as the update rule
\begin{align*}
\boldsymbol{\omega}^x_{i,j} & +=  - (\mathbf{i}_{\boldsymbol{X}} \mathbf{d} \boldsymbol{\omega})^x_{i,j} - (\mathbf{d} \mathbf{i}_{\boldsymbol{X}} \boldsymbol{\omega})^x_{i,j} \\
\boldsymbol{\omega}^y_{i,j} & +=  - (\mathbf{i}_{\boldsymbol{X}} \mathbf{d} \boldsymbol{\omega})^y_{i,j} - (\mathbf{d} \mathbf{i}_{\boldsymbol{X}} \boldsymbol{\omega})^y_{i,j}.
\end{align*}

\begin{figure}[h]\vspace*{-14mm}
\begin{picture}(200,200)(0,0)
\put(0,0){\includegraphics[width=\columnwidth]{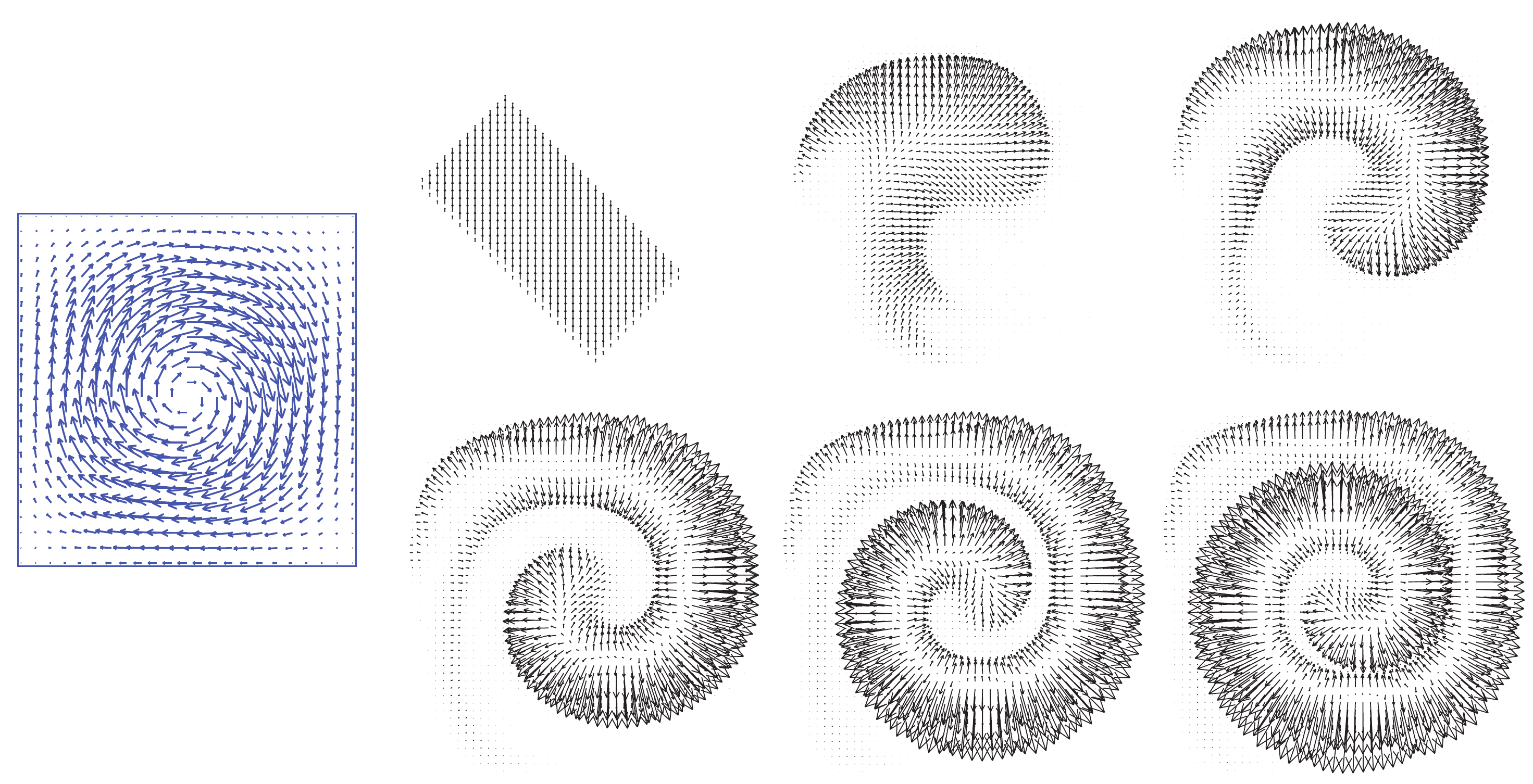}}%
\put(1,35){Rudman vortical field}
\end{picture}\vspace*{-1mm}
\caption{\textbf{High-order Advection:} in a vortical vector field (left) typically used for scalar advection, a piecewise-constant form is advected on a unit square periodic domain with a grid resolution of $48^2$ and a time step $dt\!=\!10^{-3}$ for 0, 200, and 400 steps (top), 600, 800, and 1000 steps (bottom).\vspace*{-2mm}} \label{fig:advection7}
\end{figure}

\begin{SCfigure}[1.0][h]
\begin{minipage}{0.45\textwidth}
\includegraphics[width=\textwidth]{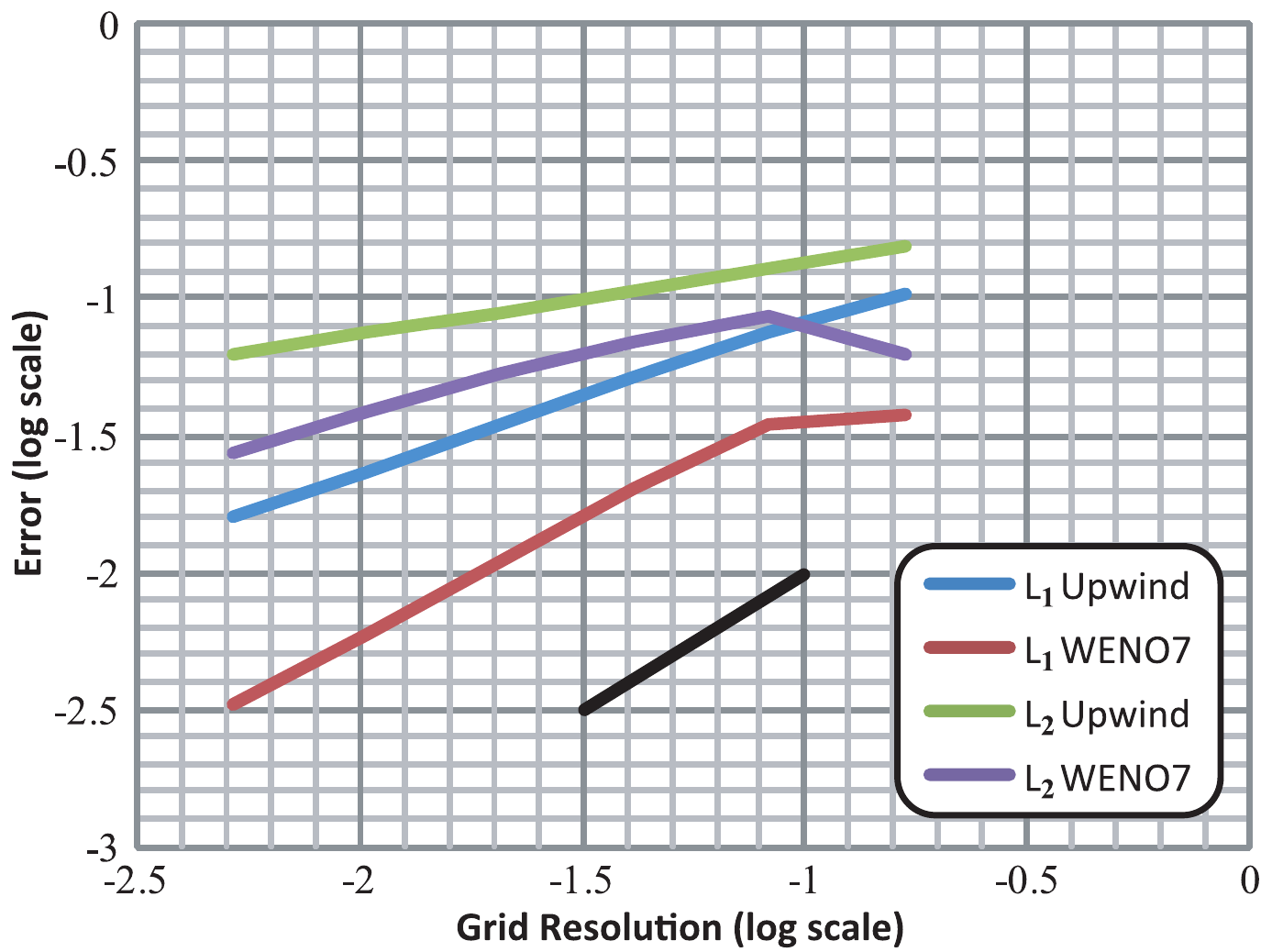}\\\vspace*{-.5mm}
\hspace*{.5\textwidth}(c)\vspace*{-5mm}
\end{minipage}
\caption{\textbf{Error Plots and Convergence Rates:} We provide error plots in $L_1$ and $L_2$ norms for power-of-two refinements
for (a) advection of a smooth form in a constant vector field, (b) advection of a smooth form in the vortical vector field of Figure~\ref{fig:advection7}(left), and (c) advection of the discontinuous form used in Figure~\ref{fig:advection7}. The black bold segments indicate a slope of $1$.} \label{fig:errorPlots}\vspace*{1mm}
\end{SCfigure}

\paragraph{A First Example.} An example of this low-order scheme can be seen in Figure~\ref{fig:compare0-7}(a-b) where we advect a piecewise constant $1$-form by a constant vector field $\boldsymbol{X} = (1,1)$ in a periodic domain.  Advecting the form forward in this velocity field for a time of $1s$ brings the form back to its original position in the continuous case; however, our numerical scheme proves very diffusive, as expected on discontinuous forms. We can however measure the error of our scheme by comparison with initial conditions as a function of the grid resolution with appropriately scaled time step sizes.  To measure the error we recall the $L_p$ norm of a $k$-form $\boldsymbol{\omega}$ is defined over a smooth manifold $\boldsymbol{\mathcal{M}}$  as $$|\boldsymbol{\omega}|_p = \left[ \int_{\boldsymbol{\mathcal{M}}} |\boldsymbol{\omega}|^p d \mu \right]^{1/p} \quad \text{where} \quad |\boldsymbol{\omega}| = (\boldsymbol{\omega},\boldsymbol{\omega})^{\frac{1}{2}}_{\boldsymbol{\mathcal{M}}}$$ and $(\cdot,\cdot)_{\boldsymbol{\mathcal{M}}}$ is the scalar product of $k$-forms defined by the Riemannian metric, and $d\mu$ is its associated volume form.  We hence define the $1$- and $2$-norms of discrete $1$-forms on a 2D regular grid with spacing $h$ as
\begin{align*}
|\boldsymbol{\omega}|_1 & = h \sum_{i,j} (|\boldsymbol{\omega}^x_{i,j}| + |\boldsymbol{\omega}^y_{i,j}|) \\
|\boldsymbol{\omega}|_2 & = \left(\sum_{i,j} (|\boldsymbol{\omega}^x_{i,j}|^2 + |\boldsymbol{\omega}^y_{i,j}|^2)\right)^{\frac{1}{2}}
\end{align*}
for simplicity, but we found using more sophisticated discretizations of the norms all yielded similar results.  Figure~\ref{fig:errorPlots}(c) shows the error plot in $L_1$ and $L_2$ norms of this simple example under power-of-two refinement, confirming the first-order accuracy of our approach.

\paragraph{High-Resolution Methods}
Note that had we chosen to leverage more sophisticated finite volume solvers in the previous example, the only changes would occur in Equations~\eqref{eq:ixdExample} and~\eqref{eq:ixExample} which would use the new numerical flux for computing the discrete contraction: any 2D method could be used for Equation~\eqref{eq:ixdExample}, while a 1D method is required for Equation~\eqref{eq:ixExample}. Due to the dimensional splitting obtaining higher order schemes is not easy, but for many application the order of accuracy is not always the most important thing.  In particular, in the presence of discontinuous solutions \emph{high-resolution} methods are often preferred for their ability to better preserve discontinuities and reduce diffusion. To test the utility and effectiveness of such schemes applied to forms, we compare the piecewise-constant upwinding method from the previous section with a Finite Volume $7^\text{th}$-order 1D WENO upwind scheme (see an overview of FV-WENO methods in~\cite{ShuOverview}).  Figure~\ref{fig:compare0-7}(c) shows the high-resolution finite volume scheme does a much better job at preserving the discontinuities, despite both methods being of the same order of accuracy for this discontinuous initial form (Figure~\ref{fig:errorPlots}(c)).

\begin{figure}[h]
\centerline{\includegraphics[width=0.472\textwidth]{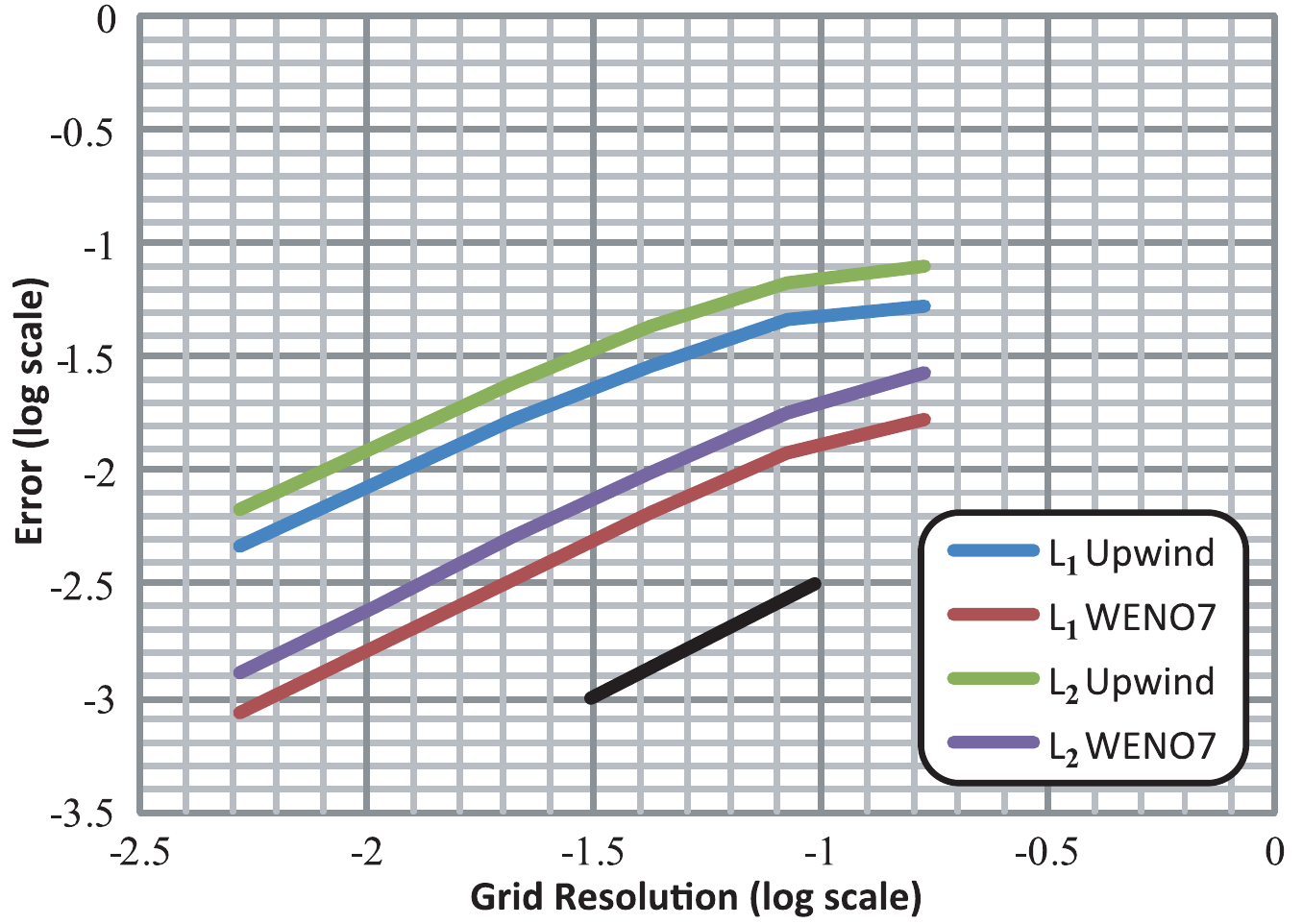}
\includegraphics[width=0.45\textwidth]{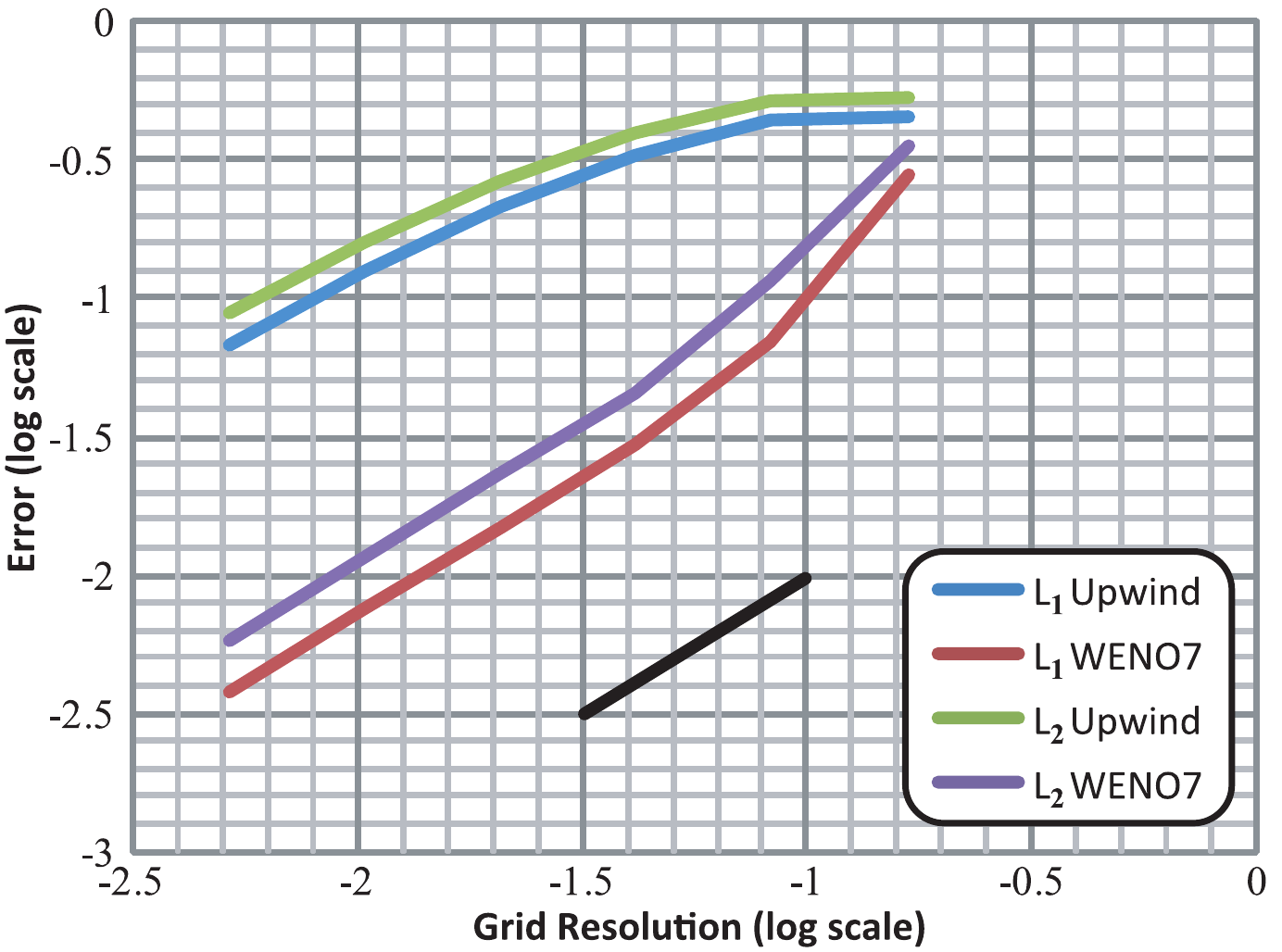}\vspace*{-.5mm}}
\hspace*{.27\textwidth}(a)\hspace*{.43\textwidth}(b)\vspace*{-1mm}
\end{figure}\vspace*{-6mm}

\paragraph{Accuracy}
To further demonstrate that properties from the underlying finite volume schemes chosen (including their accuracy) carry over to the advection of forms, we provide additional numerical tests. In Figure~\ref{fig:advection7}, we advected a simple discontinuous $1$-form in a vortical shearing vector field (Rudman vortex, left) on a $48$x$48$ grid representing a periodic domain. As expected, the form is advected in a spiral-like fashion. By advecting the shape back in the the negated velocity field for the same amount of time, we can derive error plots to compare the $L_1$ and $L_2$ norms for this example under refinement of the grid; see Figure \ref{fig:errorPlots} for other convergence tests when a first-order upwind method or a WENO-7 method is used in our numerical technique.

\subsection{Properties}
It is easy to show that our discrete Lie derivative will commute with the discrete exterior derivative as in the continuous case, by using Cartan's formula and a discrete exterior derivative which satisfies $\mathbf{d}\mathbf{d} = 0$ since we have
\begin{align*}
\mathbf{d} \boldsymbol{\mathcal{L}_X \omega} & = \mathbf{d} (\mathbf{i}_{\boldsymbol{X}} \mathbf{d} \boldsymbol{\omega} +  \mathbf{d} \mathbf{i}_{\boldsymbol{X}} \boldsymbol{\omega}) \\
& = \mathbf{d} \mathbf{i}_{\boldsymbol{X}} \mathbf{d} \boldsymbol{\omega} \\
& = (\mathbf{d} \mathbf{i}_{\boldsymbol{X}} + \mathbf{i}_{\boldsymbol{X}} \mathbf{d}) \mathbf{d} \boldsymbol{\omega} \\
& = \boldsymbol{\mathcal{L}_X} \mathbf{d} \boldsymbol{\omega}.
\end{align*}
This commutativity does not depend on any properties of the discrete contraction and therefore holds regardless of the underlying finite volume scheme chosen.  A useful consequence of this fact is that \emph{the discrete Lie advection of closed forms will remain closed by construction}; \emph{i.e.}, the advection of a gradient (resp., curl) field will remain a gradient (resp., curl) field.
\medskip

Unfortunately other properties of the Lie derivative do not carry over to the discrete picture as easily.  The product rule for wedge products, for example, does not hold for the discrete wedge products defined in~\cite{Hirani:2003}, although perhaps a different discretization of the wedge product may prove otherwise.  However, the nonlinearity of the discrete contraction operator along with the upwinding potentially picking different directions on simplices and their subsimplices makes designing a discrete analog satisfying this continuous property challenging.

\section{Conclusions}
\label{sec:conclusion}

In this paper we have introduced an extension of classical finite volume techniques for hyperbolic conservation laws to handle arbitrary discrete forms. A class of first-order finite-volume-based discretizations of both contraction and Lie derivative of arbitrary forms was presented, extending Discrete Exterior Calculus to include approximations to these operators. Low numerical diffusion is attainable through the use of high-resolution finite volume methods.  The advection of forms and vector fields are applicable in a multitude of problems, including conservative interface advection and conservative vorticity evolution.

\smallskip
Although finite volume methods can offer high resolution at a relatively low computational cost, numerical diffusion is still present and can accumulate over time. In addition, the numerical scheme we presented is not variational in nature, \emph{i.e.}, it is not (a priori) derived from a variational principle. These limitations are good motivations for future work.

\smallskip
While we have given numerical evidence demonstrating the manner in which resolution and accuracy is inherited from the underlying finite volume scheme, a formal analysis of stability and convergence remains to be performed.  In particular, it is desirable to understand what the stability of the underlying one-dimensional scheme implies about the resulting stability of our method.  Although norms for discrete differential forms have been defined, such tools are not always suitable for the analysis of nonlinear methods in multiple dimensions, even for scalar advection.  Such an investigation would thus be an important next step for the present work.

\smallskip
In the future, we also expect that extensions can be made to make truly high-order and high-resolution discretizations of the contraction and Lie derivative through $n$-dimensional reconstructions of $k$-forms and extrusions.  In particular, this would greatly facilitate the extension to simplicial meshes.  While there has been recent progress on high-order schemes for triangular meshes~\cite{TriangleWENO,TriangleCentralWENO}, these are not directly applicable for contractions of arbitrary forms.  Hence, despite having a well-defined discrete exterior derivative for simplicial meshes, such an extension will require more than just $0$- or $n$-form advection schemes, as our current dimension splitting approach to generalize such schemes for contraction of forms of other degrees does not immediately apply on non-rectangular meshes.

\begin{acknowledgement}
This research was partially supported by NSF grants CCF-0811313 \& 0811373 \& 1011944, CMMI-0757106 \& 0757123 \& 0757092, IIS-0953096, and DMS-0453145, and by the Center for the Mathematics of Information at Caltech.
\end{acknowledgement}

\bibliographystyle{plain}
\bibliography{hola}

\begin{thebibliography}{10}

\bibitem{AMR}
R.~Abraham, J.~E. Marsden, and T.~Ratiu.
\newblock {\em {Manifolds, Tensor Analysis, and Applications}}.
\newblock Applied Mathematical Sciences Vol. 75, Springer, 1988.

\bibitem{IMA_book}
D.~N. Arnold, P.~B. Bochev, R.~B. Lehoucq, R.~A. Nicolaides, and M.~Shashkov,
  editors.
\newblock {\em Compatible Spatial Discretizations}, volume 142 of {\em I.M.A.
  Volumes}.
\newblock Springer, 2006.

\bibitem{Arnold:2006:FEEC}
D.~N. Arnold, R.~S. Falk, and R.~Winther.
\newblock Finite element exterior calculus, homological techniques, and
  applications.
\newblock {\em Acta Numerica}, 15:1--155, 2006.

\bibitem{ArFaWi2010}
D.~N. Arnold, R.~S. Falk, and R.~Winther.
\newblock Finite element exterior calculus: from {H}odge theory to numerical
  stability.
\newblock To appear in \emph{Bull.\ Amer.\ Math.\ Soc.}, 74 pages, 2010.

\bibitem{Bochev2006}
P.~B. Bochev and J.~M. Hyman.
\newblock {Principles of Mimetic Discretizations of Differential Operators}.
\newblock {\em I.M.A. Volumes}, 142:89--119, 2006.

\bibitem{BossavitBook}
A.~Bossavit.
\newblock {\em Computational Electromagnetism}.
\newblock Academic Press, Boston, 1998.

\bibitem{Bossavit:2003}
A.~Bossavit.
\newblock {Extrusion, Contraction : their Discretization via Whitney Forms.}
\newblock {\em COMPEL: The International Journal for Computation and
  Mathematics in Electrical and Electronic Engineering}, 22(3):470--480, 2003.

\bibitem{Burke}
W.~L. Burke.
\newblock {\em Applied Differential Geometry}.
\newblock Cambridge University Press, 1985.

\bibitem{Carroll}
S.~Carroll.
\newblock {\em Spacetime and Geometry: An Introduction to General Relativity}.
\newblock Pearson Education, 2003.

\bibitem{Cartan1945}
{\'E}.~Cartan.
\newblock {\em Les Syst\`emes Differentiels Exterieurs et leurs Applications
  G\'eometriques}.
\newblock Hermann, Paris, 1945.

\bibitem{DKT06}
M.~Desbrun, E.~Kanso, and Y.~Tong.
\newblock {Discrete Differential Forms for Computational Sciences}.
\newblock In Eitan Grinspun, Peter Schr{\"o}der, and Mathieu Desbrun, editors,
  {\em Discrete Differential Geometry}, Course Notes. ACM SIGGRAPH, 2006.

\bibitem{Dupont2003}
T.~F. Dupont and Y.~Liu.
\newblock {Back-and-Forth Error Compensation and Correction Methods for
  Removing Errors Induced by Uneven Gradients of the Level Set Function}.
\newblock {\em {Journal of Computational Physics}}, 190(1):311–--324, 2003.

\bibitem{MOF}
V.~Dyadechko and M.~Shashkov.
\newblock {Moment-of-Fluid Interface Reconstruction}.
\newblock LANL Technical Report LA-UR-05-7571, 2006.

\bibitem{ETKSD07}
S.~Elcott, Y.~Tong, E.~Kanso, P.~Schr\"{o}der, and M.~Desbrun.
\newblock Stable, circulation-preserving, simplicial fluids.
\newblock {\em ACM Trans. Graph.}, 26(1):4, 2007.

\bibitem{EngquistOsher}
B.~Engquist and S.~Osher.
\newblock One-sided difference schemes and transonic flow.
\newblock {\em PNAS}, 77(6):3071--3074, 1980.

\bibitem{Flanders}
H.~Flanders.
\newblock {\em Differential Forms and Applications to Physical Sciences}.
\newblock Dover Publications, 1990.

\bibitem{Frankel}
T.~Frankel.
\newblock {\em The Geometry of Physics}.
\newblock Second Edition. Cambridge University Press, United Kingdom, 2004.

\bibitem{GY03}
X.~Gu and S.-T. Yau.
\newblock {Global Conformal Surface Parameterization}.
\newblock In {\em Symposium Geometry Processing}, pages 127--137, 2003.

\bibitem{GeoIntBook}
E.~Hairer, C.~Lubich, and G.~Wanner.
\newblock {\em Geometric Numerical Integration: Structure-Preserving Algorithms
  for ODEs}.
\newblock Springer, 2002.

\bibitem{MACgrids}
F.~H. Harlow and J.~E. Welch.
\newblock {Numerical Calculation of Time-dependent Viscous Incompressible Flow
  of Fluid with Free Surfaces}.
\newblock {\em Phys. Fluids}, 8:2182–--2189, 1965.

\bibitem{Heumann08}
H.~Heumann and R.~Hiptmair.
\newblock {Extrusion contraction upwind schemes for convection-diffusion
  problems}.
\newblock {Seminar f\"ur Angewandte Mathematik} SAM 2008-30, {ETH Z{\"u}rich},
  October 2008.

\bibitem{Pullin}
D.~J. Hill and D.~I. Pullin.
\newblock {Hybrid Tuned Center-Difference-WENO Method for Large Eddy
  Simulations in the Presence of Strong Shocks}.
\newblock {\em J. Comput. Phys.}, 194(2):435--450, 2004.

\bibitem{Hiptmair02}
R.~Hiptmair.
\newblock Finite elements in computational electromagnetism.
\newblock {\em Acta Numerica}, 11:237--339, 2002.

\bibitem{Hirani:2003}
A.~N. Hirani.
\newblock {\em Discrete Exterior Calculus}.
\newblock PhD thesis, Caltech, May 2003.

\bibitem{Iske2004}
A.~Iske and M.~K{\"a}ser.
\newblock {Conservative Semi-Lagrangian Advection on Adaptive Unstructured
  Meshes}.
\newblock {\em {Numerical Methods for Partial Differential Equations}},
  20(3):388--411, 2004.

\bibitem{clawpack}
R.~J. Leveque.
\newblock {CLAWPACK, at \url{http://www.clawpack.org}}.
\newblock 1994-2009.

\bibitem{LevequeBook}
R.~J. LeVeque.
\newblock {\em {Finite Volume Methods for Hyperbolic Problems }}.
\newblock {Cambridge Texts in Applied Mathematics}. {Cambridge University
  Press}, 2002.

\bibitem{TriangleCentralWENO}
D.~Levy, S.~Nayak, C.W. Shu, and Y.~T. Zhang.
\newblock {Central WENO Schemes for Hamilton-Jacobi Equations on Triangular
  Meshes}.
\newblock {\em J. Sci. Comput.}, 27:532--552, 2005.

\bibitem{WENO}
X.~D. Liu, S.~Osher, and T.~Chan.
\newblock {Weighted Essentially Non-oscillatory Schemes}.
\newblock {\em J. Sci. Comput.}, 126:202--212, 1996.

\bibitem{LovelockRund}
D.~Lovelock and H.~Rund.
\newblock {\em Tensors, Differential Forms, and Variational Principles}.
\newblock Dover Publications, 1993.

\bibitem{WestMarsden}
J.~E. Marsden and M.~West.
\newblock {Discrete Mechanics and Variational Integrators}.
\newblock {\em Acta Numerica}, 2001.

\bibitem{MoritaBook}
S.~Morita.
\newblock {\em Geometry of Differential Forms}.
\newblock Translations of Mathematical Monographs, Vol. 201. Am. Math. Soc.,
  2001.

\bibitem{Munkres1984}
J.~R. Munkres.
\newblock {\em {Elements of Algebraic Topology}}.
\newblock Addison-Wesley, Menlo Park, CA, 1984.

\bibitem{Nedelec86}
J.-C. N{\'e}d{\'e}lec.
\newblock {Mixed Finite Elements in 3D in H(div) and H(curl)}.
\newblock {\em Springer Lectures Notes in Mathematics}, 1192, 1986.

\bibitem{Nicolaides1997}
R.~A. Nicolaides and X.~Wu.
\newblock {Covolume Solutions of Three Dimensional Div-Curl Equations}.
\newblock {\em SIAM J. Numer. Anal.}, 34:2195, 1997.

\bibitem{OsherFedkiwBook}
S.~Osher and R.~Fedkiw.
\newblock {\em Level Set Methods and Dynamic Implicit Surfaces}, volume 153 of
  {\em Applied Mathematical Sciences}.
\newblock Springer-Verlag, New York, 2003.

\bibitem{SethianBook}
J.~A. Sethian.
\newblock {\em Level Set Methods and Fast Marching Methods}, volume~3 of {\em
  Monographs on Appl. Comput. Math.}
\newblock Cambridge University Press, Cambridge, 2nd edition, 1999.

\bibitem{Shi:02:FVWENO}
J.~Shi, C.~Hu, and C.~W. Shu.
\newblock {A technique for treating negative weights in WENO schemes}.
\newblock {\em J. Comput. Phys.}, 175:108--127, 2002.

\bibitem{ShuOverview}
C.-W. Shu.
\newblock {\em Essentially non-oscillatory and weighted essentially
  non-oscillatory schemes for hyperbolic conservation laws}, volume 1697 of
  {\em Lecture Notes in Mathematics}, pages 325--432.
\newblock Springer, 1998.

\bibitem{ENO}
C.~W. Shu and S.~Osher.
\newblock {Efficient Implementation of Essentially non-Oscillatory Shock
  Capturing Schemes}.
\newblock {\em J. Sci. Comput.}, 77:439--471, 1988.

\bibitem{StToDeMa2008}
A.~Stern, Y.~Tong, M.~Desbrun, and J.~E. Marsden.
\newblock Variational integrators for {M}axwell's equations with sources.
\newblock In {\em Progress in Electromagnetics Research Symposium}, volume~4,
  pages 711--715, June 2008.

\bibitem{Titarev:04:FVWENO}
V.~A. Titarev and E.~F. Toro.
\newblock {Finite-volume WENO schemes for three-dimensional conservation laws}.
\newblock {\em J. Comput. Phys.}, 201(1):238--260, 2004.

\bibitem{TACSD07}
Y.~Tong, P.~Alliez, D.~Cohen-Steiner, and M.~Desbrun.
\newblock {Designing Quadrangulations with Discrete Harmonic Forms}.
\newblock In {\em Proc. Symp. Geometry Processing}, pages 201--210, 2006.

\bibitem{Whitney1957}
H.~Whitney.
\newblock {\em Geometric Integration Theory}.
\newblock Princeton Press, Princeton, 1957.

\bibitem{TriangleWENO}
Y.~T. Zhang and C.~W. Shu.
\newblock {High-Order WENO Schemes for Hamilton-Jacobi Equations on Triangular
  Meshes}.
\newblock {\em J. Sci. Comput.}, 24:1005--1030, 2003.

\end{thebibliography}

\end{document}